\documentclass{compositio}
\usepackage[latin1]{inputenc}
\usepackage{amsmath}
\usepackage{amscd}
\usepackage{mathrsfs}
\usepackage[all]{xy}
\usepackage{graphicx}
\usepackage{pstricks} 
\usepackage{amsfonts}   
\usepackage{pict2e}   
\usepackage{makeidx}
\usepackage[colorlinks=true,linkcolor=red,citecolor=red]{hyperref}

\usepackage{amssymb}
\usepackage{latexsym}

\makeindex

\newcommand{\bs}[1]{\boldsymbol{#1}}

\newcommand{\e}{\bs{\mathrm{e}}}

\renewcommand{\H}{\mathscr{H}}

\newcommand{\M}{\mathrm{M}}
\renewcommand{\L}{\mathcal{L}}

\newcommand{\R}{\mathcal{R}}

\newcommand{\sm}[1]{\begin{smallmatrix}#1\end{smallmatrix}}

\newcommand{\simto}{\stackrel{\sim}{\to}}

\newtheorem{theorem}{Theorem}[section]
\newtheorem{proposition}[theorem]{Proposition}

\newtheorem{definition}[theorem]{Definition}
\newtheorem{remark}[theorem]{\textsc{Remark}}

\newtheorem{hypothesis}[theorem]{\textsc{Hypothesis}}

\numberwithin{equation}{section}

\title[Algorithm computing non solvable spectral radii]{An 
algorithm computing non solvable spectral radii of $p$-adic 
differential equations.}

\author{Andrea Pulita}
\email{pulita@math.univ-montp2.fr}
\address{Département de Mathématiques,
Université de Montpellier II, CC051,
Place Eugène Bataillon,
F-34 095, Montpellier CEDEX 5.}

\date{\today}



\begin{abstract}
Nous obtenons un algorithme pour le calcul explicite des valeurs 
des rayons de convergence spectrales non solubles des solutions 
d'un module différentiel sur un point de type 2, 3 ou 4 de la droite 
affine de Berkovich.
\end{abstract}

\begin{abstract}
We obtain an algorithm computing explicitly the values of 
the non solvable spectral radii of convergence of the solutions of 
a differential module over a point of type 2, 3 or 4 of the Berkovich 
affine line.
\end{abstract}

\begin{document}
\maketitle


\if{\makeatletter
\renewcommand\tableofcontents{%
    \subsection*{\contentsname}%
    \@starttoc{toc}%
    }
\makeatother

\begin{small}
\setcounter{tocdepth}{3} \tableofcontents
\end{small}
}\fi

\setcounter{section}{0}

\section*{\textsc{Introduction}}
\addcontentsline{toc}{section}{\textsc{Introduction}}
By a theorem of Young \cite{Young} the small values of the radii of convergence of the solutions of a differential operator are explicit,
and coincides with the small slopes of the Newton polygon of a differential operator attached to the module 
(cf. Prop.\ref{Prop : slopes M = slopes L}). 
Larger radii are not immediately readable on the coefficients of the operator. 
This discrepancy is a peculiarity of the p-adic world: only small radii are ``\emph{visible}''. 
To overcome this problem B.Dwork observed that the pull-back by Frobenius functor increase the radii of the solutions, and together with 
G.Christol \cite{CH-DW} they constructed an inverse of the Frobenius functor (often called Frobenius antecedent) 
in order to make the radii of the solutions smaller, and hence explicitly intelligible in a cyclic basis. 
Although theoretically satisfactory, the inversion of Frobenius is a completely implicit operation. 
Moreover the antecedent exists only if all the radii of the solutions are not small. So one is obliged to factorize the module by the radii of the 
solutions if one wants to understand the non minimal radii of the solutions. The factorization is also an implicit operation. 
Recently in \cite{Kedlaya-Book} K.Kedlaya observed that the Frobenius Push-forward operation has essentially 
the same effects as the inversion of the Frobenius on the radii of the solutions, and he is able to control the exact behavior of all the radii of the 
solutions under this operation (even small radii).\footnote{The very first reference for this is \cite{Christol-GEAU}, in which the author introduces the 
push-forward, and its relation with the pull-back, and use it to prove  the inversion of Frobenius (existence of the antecedent) under the condition that the 
solutions of the module are Bounded.} Frobenius push-forward functor 
is completely explicit, and it allows to obtain a concrete algorithm to compute the radii of the solutions that are not maximal (i.e. non solvable). 
The price to pay is that the dimension of the push-forward by Frobenius is $p$-times that of the original module, 
so that the complexity of the algorithm is multiplied by $p$ at each application of the push-forward.\footnote{Xavier Caruso recently pointed 
out that  the explicit factorization of an operator by the radii of convergence seems to be concretely implementable into a machine. 
This would highly reduce the complexity of the present algorithm since by considering the right factor of the push-forward by Frobenius 
the dimension remains constant at each step.}
Hence more the radius is large, more the complexity increase.
Moreover the algorithm admits an end if and only if all the radii are not maximal (i.e. non \emph{solvable}).
Eventually we provide the algorithm, but we avoid to provide a complete formula as one does for rank one equations (cf. \cite{F}\footnote{The 
formula that we have contributed to prove in \cite{F} is based on a completely different approach, and it uses Witt vectors 
(following techniques of \cite{Rk1}) to explicitly describe Taylor solutions of a rank one differential equation.}).
Indeed the complexity seems so great that the formula would result not useful to be written. 
Explicit examples are quite complicate even in the rank one case. 
This note is intended to make explicit the computations, 
and the link between the different results, 
in view to make them explicitly calculable by a computer.
\subsubsection{Acknowledgements.} We thanks  Kiran S. Kedlaya for Remark \ref{Kedlaya}.

\section{Radii of convergence}
Let $(K,|.|)$ be a complete field with respect to an 
ultrametric absolute value $|.|:K\to\mathbb{R}_{\geq 0}$.
Let $L/K$ be a complete valued field extension, let $c\in L$, and $\rho>0$. 
For all polynomial $P(T):=\sum_ia_iT^i\in K[T]$ define 
$|P|_{c,\rho}:=\sup_{n\geq 0}|P^{(i)}(c)|\rho^n/n!$. 
 The setting $|P_1/P_2|_{c,\rho}:=
|P_1|_{c,\rho}/|P_2|_{c,\rho}$ defines an absolute value 
on the field of fraction $K(T):=\mathrm{Frac}(K[T])$, and hence a Berkovich point $\xi_{c,\rho}$ of the affine line 
$\mathbb{A}^{1,\mathrm{an}}_K$. 
Since $\rho>0$ if $c\in K$ one obtains in this way all the points of type 2 or 3 of $\mathbb{A}^{1,\mathrm{an}}_K$, 
if one allows $c\notin K$ one also has all points of type 4. 
The derivative $d/dT$ extends by continuity to the completion $\H_{c,\rho}$ of $(K(T),|.|_{c,\rho})$. 
A differential module over $\H_{c,\rho}$ is a 
finite dimensional $\H_{c,\rho}$-vector space $\M$ together with a $K$-linear 
map $\nabla:\M\to\M$ satisfying $\nabla(fm)=d(f)m+f\nabla(m)$, 
for all $f\in \H_{c,\rho}$, $m\in \M$.
Let $r(\xi_{c,\rho})\geq \rho$ be the radius of the point $\xi_{c,\rho}$ (cf. \cite[Section 1.3.1]{NP-1}). 
If $\Omega / K$ is a complete valued field extension such that $\mathbb{A}^{1,\mathrm{an}}_{\Omega}$ has an $\Omega$-rational 
point $t_{c,\rho}\in\Omega$ lifting $\xi_{c,\rho}$, then $r(\xi_{c,\rho})$ is the radius of the largest open disk 
$\mathrm{D}^-_{\Omega}(t_{c,\rho},r(\xi_{c,\rho}))$ satisfying 
$\mathrm{D}^-_{\Omega}(t_{c,\rho},r(\xi_{c,\rho}))\cap K^{\mathrm{alg}}=\emptyset$. 

\begin{definition}\label{defi : R_i}
Let $r:=\mathrm{rank}(\M)$. For $i=1,\ldots,r$ we denote by $\R_i=\R_i^{\M,\mathrm{sp}}(\xi_{c,\rho}) \leq r(\xi_{c,\rho})$ the radius of the largest 
open disk in $\mathrm{D}^-_{\Omega}(t_{c,\rho},r(\xi_{c,\rho}))$ 
centered at $t_{c,\rho}$ on which $\M$ has at least $r-i+1$ $\Omega$-linearly independent Taylor solutions 
(cf. \cite[section 4.2]{NP-1} or \cite[11.9]{Kedlaya-Book}). One has $\R_1\leq \R_2\leq \cdots\leq \R_r$. 
\end{definition}

We say that $\R_i$ is \emph{solvable} if $\R_i=r(\xi_{c,\rho})$. In this paper we provide an algorithm computing non solvable radii.

\section{Comparison of Newton polygons and computation of small radii}
Let $r\geq 1$ be an integer. A \emph{slope sequence} is the data of 
$r$ real numbers $s_1\leq\ldots\leq s_r$ in increasing order. Define the $i$-th partial height as $h_i:=s_1+\cdots+s_i$.
A slope sequence defines univocally a convex function $h:[0,r]\to\mathbb{R}$ by $h(0):=0$, 
$h(i):=h_i$, and $h(x)=s_ix+(h_i-i\cdot s_i)$ for all $x\in]i-1,i]$, $i=0,\ldots,r$. 
The function $h$ is called the \emph{Newton polygon with slopes $s_1\leq\ldots\leq s_r$}.
\begin{definition}
The Newton polygon with slopes $s_i:=s_i^{\M,\mathrm{sp}}(\xi_{c,\rho}):=\ln(\R_i^{\M,\mathrm{sp}}(\xi_{c,\rho}))$ is 
called the \emph{spectral Newton polygon} of $\M$. We denote 
by $h_i:=h_i^{\M,\mathrm{sp}}(\xi_{c,\rho})$ its $i$-th partial height.
\end{definition}
Let $\L=\sum_{i=0}^{r}g_{r-i}d^i$, $g_i\in\H_{c,\rho}$, be a differential operator with $g_0=1$ and $g_r\neq 0$. 
Let $v_0=0$, and for all $i=1,\ldots,r$ let $v_i:=-\ln(|g_i|_{c,\rho}/\omega^i)$, where $\omega:=\lim_n|n!|^{1/n}$. 
Let $L_i:=\{(x,y)\in\mathbb{R}^2\;\textrm{such that }\;x=i,y\geq v_i\}$, note that $L_i$ is empty if and only if $g_i=0$. 
Define the \emph{spectral Newton polygon} $NP(\L)$ as the 
intersection of all upper half planes 
$H_{a,b}:=\{(x,y)\in\mathbb{R}^2\textrm{ such that }y\geq ax+b\}$ with  
$\{(i,v_i)\}_{i=0,\ldots,r}\subset H_{a,b}$.
Let $h^{\L}:[0,r]\to\mathbb{R}$ be the convex function whose epigraph is $NP(\L)$: 
$h^{\L}(x)=\min\{y\textrm{ such that }(x,y)\in NP(\L)\}$. 
Explicitly one has $h^\L_i:=h^{\L}(i)=\sup_{s\in\mathbb{R}}\{s\cdot i +\min_{j=0,\ldots,r}v_j-s\cdot j\}$. 
Then  $NP(\L)$ is the Newton polygon with slopes $\{s_i^{\L}:=h^{\L}_i-h^{\L}_{i-1}\}_{i=1,\ldots,r}$.
\begin{proposition}[(\protect{\cite{Young}})]\label{Prop : slopes M = slopes L}
Let $\L$ be a differential operator as above and let $\M$ be the differential module defined by $\L$. 
Let $C:=\ln(\omega\cdot r(\xi_{c,\rho}))$, then for all $i=1,\ldots,r$ one has
\begin{equation}\label{(2.1)}
\min(s_i^{\M,\mathrm{sp}},C)\;=\;\min(s_i^{\L},C)\;.\qquad\Box
\end{equation}
\end{proposition}
\begin{remark}\label{Rk p-adic or not}
In order to apply \eqref{(2.1)} we need an algorithm to find a cyclic basis of $\M$ (cf. section \ref{cyc}).
If the absolute value of $K$ is trivial on $\mathbb{Z}$ (i.e. if $|n|=1$ for all $n\in\mathbb{Z}-\{0\}$), 
then $\omega=1$, and Proposition \ref{Prop : slopes M = slopes L} allows to find all the radii $\R_i$. 
If the absolute value of $K$ is $p$-adic (i.e. if $|p|<1$), then $\omega=|p|^{\frac{1}{p-1}}<1$. 
In this case, we also need a technique (Frobenius push-forward)  
making the (non solvable) radii smaller than $\omega r(\xi_{c,\rho})$ (cf. section \ref{Frob}).
\end{remark}
\if{
\subsection{Algorithm.}
In order to apply Prop. \ref{Prop : slopes M = slopes L} one needs
to find an explicit operator $\L$ representing $\M$ in a cyclic basis, 
this is done in section \ref{cyc} using 
Katz explicit algorithm to find an 
explicit cyclic vector \cite{Katz-cyclic-vect}.
The proposition then computes the radii $R_i(\M)$ such that 
$R_i(\M)<\omega\rho$. 

If $\omega=1$ this is enough to 

In the sequel we apply Prop. \ref{Prop : slopes M = slopes L} to 
find explicitly $s^{\M}_i$ under the assumption 
$s_i^{\M}<\ln(\rho)$.  
If $ \ln(\omega\rho)\leq s_i^{\M}<\rho$, then one replace 
$\M$ by its push-forward by Frobenius $\varphi_*\M$, one 
provides the exact relation between the $h_i^\M$ and 
$h_i^{\varphi_*(\M)}$, and shows that $s_i^{\varphi_*(\M)}$ 
fulfill 
}\fi
\section{Explicit Cyclic vector}\label{cyc}
Let $(F,d)$ be a differential field and let 
$F\langle d\rangle=\oplus_{i\geq 0}F\circ d^i$ be the Weil algebra of differential operators. 
The multiplication of $F\langle d\rangle$ extends that of $F=F\circ d^0$ by the rule 
$d\circ f=f\circ d+d(f)$, for all $f\in F$. 
Finite dimensional differential modules over $F$ are exactly torsion left $F\langle d\rangle$-modules. 
The so called \emph{cyclic vector theorem} asserts that all differential module 
are not only torsion modules over $F\langle d\rangle$ : they are cyclic modules i.e. of the form 
$(\M,\nabla)=(F\langle d\rangle/F\langle d\rangle\L,d)$, for some 
$\L:=\sum_{i=0}^rg_{r-i}d^i\in F\langle d\rangle$, with $g_0=1$, $g_i\in F$. 
The image of $\{1,d,d^2,\ldots,d^{r-1}\}$ in the quotient 
form a basis of $\M$, and the action of $\nabla$ is given by the multiplication by $d$ in the quotient. 
In fact the cyclic vector theorem is 
equivalent to the existence of an element $c\in\M$, called cyclic vector, such that $\{c,\nabla(c),\nabla^2(c),\ldots,\nabla^{r-1}(c)\}$ is a 
basis of $\M$. In this case if $c_i:=\nabla^i(c)$, and if $\nabla^r(c_0)=\sum_{i=0}^{r-1}f_ic_i$, then $f_{i}=-g_{r-i}$.
The existence of such a vector is due to \cite[Ch.II, Lemme 1.3]{Deligne-Reg-Sing}. 
Subsequently N.M.Katz provided the following explicit algorithm
\begin{theorem}[(\cite{Katz-cyclic})]\label{theorem 1 of Katz}
Let $(\M,\nabla)$ be a differential module over $(F,d)$ of rank $r$, and let $\e:=\{e_0,\ldots,e_{r-1}\}\subset\M$ 
be a basis of $\M$. Let $a_0,\ldots,a_{r(r-1)}\in F$ be $r(r-1)+1$ distinct constants i.e. $d(a_i)=0$. 
Then at least one of the following elements of $\M$ is a cyclic vector:
\begin{equation}\label{cyclic vector}
c(\e,T-a_i)\;:=\;
\sum_{j=0}^{r-1}\frac{(T-a_i)^j}{j!}
\sum_{k=0}^j(-1)^k\binom{j}{k}\nabla^k(e_{j-k})\;.\qquad\Box
\end{equation}
\end{theorem}

\begin{remark}
The explicit base change matrices to the Katz's cyclic basis are quite 
complicate and have been explicitly computed in 
\cite{Pu-CV}. 
On the other hand the proof of the existence of a cyclic vector of Deligne 
\cite[Ch.II, Lemme 1.3]{Deligne-Reg-Sing} proves that the family of 
cyclic vectors in 
a given module is the complement of an hypersurface. The base change 
matrices of Katz's algorithm are quite involved and hard to find by hand 
even in small examples. It is often convenient to pick an arbitrary vector 
and test if it is cyclic. 
\end{remark}

\begin{remark}\label{Kedlaya}
One shall avoid the use of a cyclic basis using
\cite[Lemma 6.7.3, Thm.6.7.4, Conjecture 4.4.9]{Kedlaya-Book}.
\end{remark}

\section{Frobenius Push-Forward and explicit computation of 
larger radii}\label{Frob}
In this section we assume that $|p|<1$ (cf. Remark \ref{Rk p-adic or not}).
\begin{hypothesis}
$\R_i^{\M,\mathrm{sp}}$ is insensitive to scalar extensions of $K$, and by translations. So in the sequel we will assume $c=0$ 
and replace the indexation $(c,\rho)$ by $\rho$. In this case one has $r(\xi_{\rho})=\rho$. 
We then work with $|.|_\rho$, $\xi_\rho$, $\H_\rho$, $r(\xi_\rho)=\rho$ with the evident meaning of notation. 
If the reader needs to preserve the setting $(c,\rho)$, the same 
computations hold replacing the map $\varphi:T\mapsto T^p$ by 
$T\mapsto (T-c)^{p}+c$. 
Or alternatively one also can preserve $\varphi:T\to T^p$, 
and proceed as in \cite[Section 7]{NP-1} to check the behavior of the radii by 
Frobenius at points that are close enough to the principal branch $\rho\mapsto 
|.|_{0,\rho}$ (this is often necessary if one search the slopes of $\R^{\M,\mathrm{sp}}_i$ along a 
Berkovich path $\rho\mapsto |.|_{c,\rho}$, with $c\in K$ and $\rho$ close to 
$|c|$). 
\end{hypothesis}
Let $\widetilde{T},T$ be two variables, and let $\varphi:K(T)\to K(\widetilde{T})$ be the 
ring morphism sending $T$ into $\widetilde{T}^p$. This extends into a isometric inclusion $\varphi:\H_{\rho^p}\to\H_{\rho}$ of degree $p$.
One has the rule $\frac{d}{dT}(f(T))=\frac{d/d\widetilde{T}}{p\widetilde{T}^{p-1}}(f(T))$, for all $f\in K(T)$. We call 
\begin{equation}
d_{\rho^p} := \frac{d}{dT}\;,\qquad
d_\rho := \frac{d}{d\widetilde{T}}\;,\qquad
\widetilde{d}_{\rho^p} := (p\widetilde{T}^{p-1})^{-1}\frac{d}{d\widetilde{T}}\;.
\end{equation}
Let $(\widetilde{\M},\nabla)$ is a differential module over $(\H_{\rho},d_\rho)$ of rank $r$. Since $(\widetilde{d}_{\rho^p})_{|_{\H_{\rho^p}}}=d_{\rho^p}$, 
then $(\widetilde{\M},(p\widetilde{T}^{p-1})^{-1}\nabla)$ 
is a differential module over $(\H_{\rho},\widetilde{d}_{\rho^p})$ that can be seen (by restriction of the scalars) as a differential module over 
$(\H_{\rho^p},d_{\rho^p})$ of rank $pr$. 
We call $(\varphi_*\widetilde{\M},\varphi_*\nabla)$ the differential module so obtained.

\subsection{Explicit matrix of $\varphi_*\nabla$.}
One has a direct sum decomposition $\H_\rho=\oplus_{k=0}^{p-1}\varphi(\H_{\rho^p})\cdot\widetilde{T}^k$, so that 
each $g(\widetilde{T})\in\H_\rho$ can be uniquely written as 
$g(\widetilde{T})=\sum_{k=0}^{p-1}g_k(\widetilde{T}^p)\widetilde{T}^k=\sum_{k=0}^{p-1}g_k(T)\widetilde{T}^k$.
The derivation $\widetilde{d}_{\rho^p}$ stabilizes globally each factor and 
$\widetilde{d}_{\rho^p}(g_k(T)\widetilde{T}^k)=d_{\rho^p}(g_k(T))\widetilde{T}^k$. 
For all $g(\widetilde{T})\in\H_\rho$ we define $\varphi_*(g)(T)\in M_{p\times p}(\H_{\rho^p})$ as the matrix 
of the multiplication by $g(\widetilde{T})/(p\widetilde{T}^{p-1})$, 
with respect to the basis $1,\widetilde{T},\ldots,\widetilde{T}^{p-1}$ over $\H_{\rho^p}$. One has
\begin{equation}\label{eq : phi_*(G)}
\varphi_*(g)(T)\;=\;(pT)^{-1}\cdot\left(
\sm{g_{p-1}(T)&Tg_{p-2}(T)&Tg_{p-3}(T)&\cdots&\cdots&\cdots&Tg_0(T)\\
g_{0}(T)&g_{p-1}(T)&Tg_{p-2}(T)&Tg_{p-3}(T)&\cdots&\cdots&Tg_1(T)\\
g_{1}(T)&g_{0}(T)&g_{p-1}(T)&Tg_{p-2}(T)&Tg_{p-3}(T)&\cdots&Tg_2(T)\\
\cdots&\cdots&\cdots&\cdots&\cdots&\cdots&\cdots\\
\cdots&\cdots&\cdots&\cdots&\cdots&\cdots&\cdots\\
g_{p-2}(T)&g_{p-3}(T)&\cdots&\cdots&g_1(T)&g_{0}(T)&g_{p-1}(T)}
\right)
\end{equation}
Notice that the terms over the diagonal are multiplied by $T$. 
Let $(\widetilde{\M},\widetilde{\nabla})$ be a differential module over $\H_\rho$. Fix a 
$\H_\rho$-linear isomorphism $\H_\rho^r\simto\widetilde{\M}$ (i.e. a basis of $\widetilde{\M}$), 
and let $\frac{d}{d\widetilde{T}}-G(\widetilde{T})$ be the map $\widetilde{\nabla}$ in this basis, 
where $G(\widetilde{T})=
(g_{i,j}(\widetilde{T}))_{i,j=1,\ldots,r}\in M_{r\times r}(\H_\rho)$. Writing 
$\H_\rho^r=(\oplus_{k=0}^{p-1}\varphi(\H_{\rho^p})\cdot\widetilde{T}^k)^r$ one sees that the matrix of 
$\varphi_*(\widetilde{\nabla})$ is given by the block matrix 
\begin{equation}
\varphi_*(G)(T)\;:=\; 
(\varphi_*(g_{i,j})(T))_{i,j=1,\ldots,r}\;\in\; M_{pr\times pr}(\H_{\rho^p})\;.
\end{equation}
\subsection{Behavior of the radii by Frobenius push-forward}
\begin{theorem}[(\protect{\cite[Thm.10.5.1]{Kedlaya-Book}})]\label{Kedlpush-f}
Let $\R_1\leq\ldots\leq \R_r$ be the radii of $\widetilde{\M}$ at $\xi_{\rho}$ (cf. Def.\ref{defi : R_i}). 
Let $i_1$ be such that $\R_{i_1}\leq \omega \rho<\R_{i_1+1}$.\footnote{It is understood that $i_1=0$ if $\omega\rho<\R_1$.} 
Then, up to permutation, the spectral radii of $\varphi_*\widetilde{\M}$ at $\xi_{\rho^p}$ are 
\begin{equation}
\bigcup_{i\leq i_1}\Bigl\{\underbrace{|p|\rho^{p-1}\R_i,
\ldots,|p|\rho^{p-1}\R_i}_{p\textrm{-times}}\Bigr\}
\bigcup_{i> i_1}\Bigl\{\R_i^p,
\underbrace{\omega^p\rho^p,\ldots,\omega^p\rho^p}_{p-1\textrm{-times}}\Bigr\}\;.\qquad\Box
\end{equation}
\end{theorem}
If $s_1\leq\ldots\leq s_r$ is the slope sequence of the spectral Newton polygon of 
$\widetilde{\M}$ at $\xi_{\rho}$, and if $i_0\geq i_1$ satisfies $\R_{i_0}<\rho=\R_{i_0+1}$,\footnote{It is 
understood that $i_0=r$ if $\R_r<\rho$.} then by Theorem \eqref{Kedlpush-f} the slope sequence 
associated to $\varphi_*\widetilde{\M}$ at $\xi_{\rho^p}$ is 
\begin{eqnarray}
\overbrace{\ln\Bigl(|p|\rho^{p-1}\Bigr)+s_{1}
=\cdots=\ln\Bigl(|p|\rho^{p-1}\Bigr)+s_{1}}^{p\textrm{-times}}
\leq&\cdots&\leq
\overbrace{\ln\Bigl(|p|\rho^{p-1}\Bigr)+s_{i_1}=\cdots=
\ln\Bigl(|p|\rho^{p-1}\Bigr)+s_{i_1}}^{p\textrm{-times}}\leq\qquad\nonumber\\
\leq\overbrace{\ln\Bigl(\omega^p\rho^p\Bigr)=\cdots=\ln\Bigl(\omega^p\rho^p\Bigr)}^{(p-1)(r-i_1)\textrm{-times}}
&<&p s_{i_1+1}\leq\cdots\leq p s_{i_0}<\overbrace{\ln(\rho^p)=\cdots=\ln(\rho^p)}^{(r-i_0)\textrm{-times}}.\label{(4.3)}
\end{eqnarray}
We have two main goals here. Firstly the sequence 
$s_1\leq\cdots\leq s_r$ is perfectly determined by the knowledge of the 
slope sequence \eqref{(4.3)} of 
$\varphi_*\widetilde{\M}$ (even if some of the $s_i$ are equal to the 
critical value $\ln(\omega\rho)$), see \cite[Prop. 6.17]{NP-1} for a more 
precise statement. 
Secondly the values of $s_i$ satisfying $ \ln(\omega^{1/p}\rho)\leq s_i < \ln(\omega^{1/p}\rho) $ 
corresponds to small radii\footnote{i.e. radii that are smaller 
than $\omega\rho^p$.} of $\varphi_*\widetilde{\M}$ that are explicitly 
intelligible by Proposition \ref{Prop : slopes M = slopes L}. 
Iterating this construction by performing several times the push-forward one obtains an explicit algorithm that 
computes all the non solvable radii $\R_1,\ldots,\R_{i_0}$ in a finite number of steps. 
Once this have been achieved, one knows in fact all the spectral radii since the remaining radii are all equal to $\rho$. 
Unfortunately Proposition \ref{Prop : slopes M = slopes L} does not furnish any information about radii that are larger 
that $\omega\rho$, so (unless the radii are all not solvable) it seems impossible to know whether the algorithm is ended or if one needs more 
applications of the Frobeinius push-forward.

\bibliographystyle{amsalpha}
\bibliography{bib}

\end{document}